\documentclass[12pt]{article}
\newcommand{\qed}{\ifmmode$\Box$\else{\unskip\nobreak\hfil
\penalty50\hskip1em\null\nobreak\hfil$\Box$
\parfillskip=0pt\finalhyphendemerits=0\endgraf}\fi}
\usepackage{amsmath}
\usepackage{latexsym}
\usepackage{amsfonts}
\begin{document}
\renewcommand{\theequation}{\thesection.\arabic{equation}}
\newtheorem{Pa}{Paper}[section]
\newtheorem{Tm}[Pa]{{\bf Theorem}}
\newtheorem{La}[Pa]{{\bf Lemma}}
\newtheorem{Cy}[Pa]{{\bf Corollary}}
\newtheorem{Rk}[Pa]{{\bf Remark}}
\newtheorem{Pn}[Pa]{{\bf Proposition}}
\newtheorem{Ee}[Pa]{{\bf Example}}
\newtheorem{Dn}[Pa]{{\bf Definition}}
\newtheorem{I}[Pa]{{\bf}}
\oddsidemargin 0in
\topmargin -0.5in
\textwidth 16.5truecm
\textheight 23truecm
\parindent =0.cm
\textwidth 16.5truecm
\textheight 23truecm
\begin{center}
{\bf \Large On the structure of Banach algebras associated with automorphisms}\\
\bigskip
{\large A. Lebedev  }\\
 {Institute of Theoretical Physics\\
University in Bia{\l}ystok
\\Lipowa 41, 15-424 Bia{\l}ystok, Poland}
\end{center}
\bigskip
\begin{abstract}
In the present paper we study the structure of a Banach algebra $B(A, T_g )$ generated
by a certain Banach algebra $A$ of operators acting in a Banach space $D$ and a group
$\{ T_g \}_{g \in G}$ of isometries of $D$  such that $T_g A T^{-1}_g = A$.
We investigate the interrelations
between the existence of the  expectation of $B(A, T_g )$ onto $A$, topological freedom of
the automorphisms of $A$ induced by $T_g$ and  the dual action of the group $G$ on $B(A, T_g )$.
The  results obtained are applied to the description of the structure of
Banach algebras generated
by 'weighted composition operators' acting in various spaces.
\end{abstract}

{\bf AMS Subject Classification:} 47D30, 16W20, 46H15, 46H20

\medskip
{\bf Key Words:} Banach algebras, isometries, automorphisms, topologically free action,
dual action, Banach algebras generated by weighted composition operators

\tableofcontents
\section{Introduction}
\setcounter{equation}{0} The principal  object under consideration
in this paper  is a Banach algebra $B(A, T_g )$ generated by a
certain Banach algebra $A$ of operators acting in a Banach space
$D$ and a group $\{ T_g \}_{g \in G}$ of isometries of $D$ (a
representation $g \to T_g$ of a discrete group $G$) such that
$$
T_g A T^{-1}_g = A, \hspace{10mm} g\in G
$$
which means that $T_g$ generates the automorphism ${\hat T}_g$ of $A$ given by
$$
{\hat T}_g (a) = T_g a T^{-1}_g , \hspace{10mm} a\in A.
$$
The purpose of the paper is a description  of the structure of such  algebras.\\

In the Hilbert space situation (that is in the    $C^\ast -$algebra theory) the analogous objects
 are closely related to the crossed products
(see, for example \cite{Ped}) and  description of their structure
is the theme of numerous investigations. In particular, Landstad
\cite{Land} presented the necessary and sufficient conditions (in
terms of {\em duality theory}) for a $C^\ast -$algebra to be
isomorphic to a crossed product (of an
 algebra and a locally compact group of automorphisms).
In the case of a discrete group in \cite{AnLeb},
Chapter 2 there were found the conditions for a $C^\ast -$algebra to be isomorphic to
a crossed product in terms of the group
action (the so called {\em topologically free action} (see \ref{1.6} of this paper))
and also in
terms of satisfaction of a certain  inequality ({\em property (*)} (\ref{e1.3}) of
the present  paper)
guaranteeing the existence of the  expectation  of the algebra $B(A, T_g )$ onto the
algebra $A$ (see (\ref{e1.4}), (\ref{e1.5})).\\

One of the main aims of the paper is to investigate to what extent the mentioned properties
(topologically free action, property (*) and dual action of the group) give us a possibility to
obtain the description of the structure of $B(A, T_g )$ up to isomorphism in the general
Banach space situation.\\

In Sections 2 and 3 we investigate the interrelations between these properties and also study
the naturally arising {\em property (**)} \ref{1.4}
 (the restoration of an element of $B(A, T_g )$ by its 'Fourier coefficients').
  In particular it is shown
that in many cases the properties considered are equivalent.\\

Since in the general Banach space situation we do not have a
universal object like the crossed product in the Hilbert space situation to
describe the structure of
$B(A, T_g )$ we have to specify the algebra $A$ and the isometries $T_g$.
Sections 4-8 are devoted
to the applications of the results obtained in Sections 2, 3 to the description
of the structure of
concrete Banach algebras associated with automorphisms, namely the algebras generated by
'weighted composition operators' acting  as in the spaces of continuous functions so also in
Lebesgue spaces.\\

We establish a number of isomorphism results for the algebras investigated and in addition
find out that the arising algebras are in a way qualitatively different. In particular when
considering the operators in $C(X, E), L^\infty _\mu (X, E)$ and in $ L^1 _\mu (X, E)$
 we {\em can}
calculate their norms  (Theorems \ref{2.2}, \ref{2.17} and \ref{2.26}) while for
 the operators in
$L^p _\mu (X, E), \hspace{2mm} 1<p < \infty$ we have nothing like this.
 Moreover  to  obtain the  isomorphism
theorems  for the algebras $B(A, T_g )$ in   $C(X, E), L^\infty _\mu (X, E)$, $ L^1 _\mu (X, E)$
we
do not need any information on the structure of the group of operators generating automorphisms
while this structure (namely the amenability of the group $G$) is vital when we are
investigating the operators
in $L^p _\mu (X, E), \hspace{2mm} 1<p < \infty $ (see Theorem \ref{2.10} and Remark \ref{2.33}).\\

Throughout the paper we shall use the following notation.\\
For any topological space $X$ and  normed space $B$ we  denote by $C(X, B)$
the normed space consisting  of all continuous functions on $X$ taking values in $B$ with sup norm.\\
As $X$ we shall use only completely regular spaces.\\
$D, E$ and $F$ will denote some Banach spaces.\\
$L(E)$ will denote  the Banach algebra of all linear bounded operators acting in $E$.

\section{ $B(A, T_g )$. Properties (*) and (**). Topologically free action}
\setcounter{equation}{0}

In this section we start to investigate some important properties of the algebras associated
with automorphisms.
\begin{I}
\label{1.2}
\em
Let  $B(A, T_g )$ be a Banach algebra  generated
by a certain Banach algebra $A$ of operators acting in a Banach space $D$ and a group
$\{ T_g \}_{g \in G}$ of isometries of $D$ (a representation $g \to T_g$ of
a discrete group $G$) such that
\begin{equation}
\label{e1.1}
T_g A T^{-1}_g = A, \hspace{10mm} g\in G
\end{equation}
which means that $T_g$ generates the automorphism ${\hat T}_g$ of $A$ given by
\begin{equation}
\label{e1.2}
{\hat T}_g (a) = T_g a T^{-1}_g , \hspace{10mm} a\in A.
\end{equation}
\end{I}
The concrete types of algebras $A,$ spaces $D$ and operators
$T_g$ for which the theory will be developed will be specified below.
\begin{I}
\label{1.3}
{\bf Property (*).}
\em
If $A$ is a $C^\ast -$algebra of operators containing the identity and
acting in a Hilbert space $H$ and $\{ T_g \}_{g \in G}$
is a unitary representation of a group $G$ in $H$ then the $C^\ast -$algebra generated
by  $A$ and  $\{ T_g \}_{g \in G}$ will be denoted by $C^\ast (A, T_g )$.

 When considering the $C^\ast -$algebra case we found out (see \cite{AnLeb}, Chapter 2)
  that one of
 the most important properties of the algebra $C^\ast (A, T_g )$ in the
presence of which one can obtain the deep and fruitful theory of the subject is the
next {\em property (*)}:\\

for any finite sum $b = \sum a_g T_g , \hspace{3mm}a_g \in A$ the following inequality holds
\begin{equation}
\label{e1.3}
\Vert b \Vert = \Vert  \sum a_g T_g \Vert \ge \Vert a_e \Vert ,
\end{equation}
where $e$ is the identity of the group $G$.
\end{I}
One of the aims of this paper is to clarify  the meaning  of the property (*) in
the Banach algebra situation that is to investigate the role of this property for $B(A, T_g ).$

If an algebra $B(A, T_g )$ possesses the property (*) then for every $g_0 \in G$
there is correctly defined
the mapping
\begin{equation}
\label{e1.4}
N_{g_0} : \sum a_g T_g \to a_{g_0}
\end{equation}
which can be extended up to the mapping
\begin{equation}
\label{e1.5}
N_{g_0} : B(A, T_g ) \to A.
\end{equation}
 \begin{I}
\label{1.4}
{\bf Property (**).}
\em
One more important property is the following.\\
We shall say that an algebra $B(A, T_g )$ which possesses the property (*) also possesses  the
{\em property (**)}
if
\begin{equation}
\label{e1.6}
B(A, T_g ) \ni b = 0 \hspace{2mm}{\rm iff}  \hspace{2mm}N_g (b) = 0 \hspace{2mm}
{\mbox {\rm for
 every}}  \hspace{2mm}g \in G
\end{equation}
where $N_g$ is the mapping introduced above.
\end{I}
\begin{Rk}
\label{Rk1.4}
\em
In the $C^\ast -$algebra situation we have (see \cite{AnLeb}, Theorems 12.8 and 12.4):\\

{\em if $G$ is a discrete amenable group and $C^\ast (A, T_g )$
possesses the property (*) then $C^\ast (A, T_g )$ possesses the
property (**) as well. }
\end{Rk}
In fact the presence of the properties (*) and (**) makes it possible to 'reestablish' an
element $b \in B(A, T_g )$ via its 'Fourier' coefficients $N_g (b), \hspace{2mm}g\in G$
and we shall find out further that in many reasonable situations this 'reestablishing'
can be carried out successfully.\\

Recall that in the $C^\ast -$algebra situation there is a close relation between the property (*)
and the so-called {\em topological freedom} of the action of the group of automorphisms
$\{ {\hat T}_g \}$ (see \cite{AnLeb}, 12.13 and Theorem 12.14).
This fact suggests to consider the interrelation between
the corresponding properties in the algebra $B(A, T_g )$.
\begin{I}
\label{1.5a}
\em
To start with we make a number of remarks clarifying the further choice of
objects and definitions.\\
Observe first that if
 $\{ T_g \}_{g\in G}$ is a group of isometries satisfying (\ref{e1.1}) then evidently
\begin{equation}
\label{center}
T_g {\mathbb Z}(A) T_g^{-1} = {\mathbb Z}(A)
\end{equation}
 where ${\mathbb Z}(A)$ is the center of $A$.\\
 Let $A$ be a certain Banach algebra isomorphic to $C(X, B)$ where $X$ is a
completely regular space and $B$ is a Banach algebra then
\begin{equation}
\label{zenter1}
{\mathbb Z}(A) = C(X, {\mathbb Z}(B)).
\end{equation}
Indeed, the inclusion ${\mathbb Z}(A)\supset  C(X, {\mathbb Z}(B))$ is obvious and more over
if ${\bar a}\in {\mathbb Z}(A)$ then for any $a\in A$ and every $x_0 \in X$ we have
${\bar a}(x_0 ) a(x_0 ) = a(x_0 ) {\bar a}(x_0 )$. So taking $a(x) \equiv b\in B$ we conclude
that ${\bar a}(x_0 ) \in  {\mathbb Z}(B)$ that is ${\mathbb Z}(A)\subset  C(X, {\mathbb Z}(B))$.\\
Let us assume for a while that
$$
{\mathbb Z}(B) \cong C(M)
$$
(this assumption is of course rather restrictive but it will enable us to clarify a certain idea).
Then
\begin{equation}
\label{center3}
{\mathbb Z}(A) = C(X, {\mathbb Z}(B)) \cong C(X, C(M)) \cong C(X\times M).
\end{equation}
If $X$ and $M$ are compact spaces (in which case $X\times M$ is a compact space as well)
then (\ref{center3}) along with (\ref{center}) means that the automorphisms  ${\hat T}_g$
 (\ref{e1.2})
generate the  homeomorphisms $t_g : X\times M \to X\times M$ according to the formula
\begin{equation}
\label{center4}
[{\hat T}_g (z)](y) = z(t_g^{-1} (y)), \hspace{5mm}z\in {\mathbb Z}(A), \hspace{2mm} y\in
  X\times M.
\end{equation}
Clearly if $M = \{ \cdot \}$ is a single point (that is ${\mathbb Z}(B) =  \{ cI \}$
where $I$ is the identity of B)
 then $X\times M \cong X $ and $t_g$ are homeomorphisms of $X$.
\end{I}
\begin{I}
\label{1.5}
\em
 Having in mind the reasoning presented above in what follows we shall confine ourselves
to the case
\begin{equation}
\label{center5}
{\mathbb Z}(B) =  \{ cI \}
\end{equation}
 that is we would like to work with  the 'initial' space
$X$ and not to involve additional objects as $X\times M$ and the like.\\
Obviously  if $B = L(E)$ then  (\ref{center5}) is satisfied.\\
So let $A\subset L(D)$ be a Banach algebra of operators isomorphic
to $C(X, L(E))$ where $X$ is a certain completely regular space
and $E$ and $D$ are Banach spaces (thus ${\mathbb Z}(A) \cong
C(X)$). Let $\{ T_g \}_{g\in G}$ be a group of isometries
satisfying (\ref{e1.1}).  According to (\ref{center}) the
automorphisms ${\hat T}_g$  (\ref{e1.2}) preserve the center and
henceforth we assume that their action on the center is given by
\begin{equation}
\label{e1.7}
[{\hat T}_g (z)](x) = z(t_g^{-1} (x)), \hspace{5mm}z\in {\mathbb Z}(A), \hspace{2mm} x\in   X.
\end{equation}
where $t_g : X \to X$ are some homeomorphisms of $ X $.\\
{\bf Remark.} Since   ${\mathbb Z}(A) \cong C(X)\cong C({\bar X})$ where $\bar X$ is
the Stone-Cech
compactification of $X$ the automorphisms  ${\hat T}_g$ generate some homeomorphisms
$\tau_g : {\bar X}\to {\bar X}$ given by the equations
$$
[{\hat T}_g (z)](x) = z(\tau_g^{-1} (x)), \hspace{5mm}x\in {\bar X}.
$$
But in general the homeomorphisms $\tau_g$ {\em do not preserve} the subset
$X\subset {\bar X}$.\\
{\em Example.} Let $Y$ be a {\em non compact} completely regular space and
$\bar Y$ be its Stone-Cech compactification. Consider the space $X = Y\coprod {\bar Y}$
where $\coprod$ denotes  the disjoint union of spaces. Clearly
$C(X) \cong C({\bar X})= C({\bar Y} \coprod  {\bar Y})$.
There is a one to one correspondence between
the elements of $A = C(X)$ and the pairs $(a_1 , a_2 )$ where $a_i \in C({\bar Y}),
 \hspace{2mm} i=1,2$.
Let $\hat T$ be the automorphism of $A = {\mathbb Z}(A)$ given by
${\hat T}(a_1 , a_2 ) = (a_2 , a_1 )$. Evidently this automorphism generates the homeomorphism
$\tau$ of $\bar X$ which does not preserve $X$.\\
 So   (\ref{e1.7})  means that we confine ourselves to the
situation when $\tau_g$ {\em  preserve} $X$ (and in this case their restrictions onto $X$
coincide with  $t_g$).\\
 We make this assumption simply  not to overload the exposition and not to involve
$\bar X$ instead of $X$.
\end{I}
 \begin{I}
\label{1.6}
{\bf Topologically free action.}
\em
Denote by $X_g , \hspace{2mm} g\in G  $ the set
\begin{equation}
\label{e1.8}
X_g = \{ x \in X : t_g (x) = x \}.
\end{equation}
We say that the group $G$ acts {\em topologically freely} on $A$
by automorphisms ${\hat T}_g$ (or on $X$ by homeomorphisms $t_g$
defined in \ref{1.5}) if for any $g\in G,\hspace{2mm} g\ne e$ the
set $ X_{g} $
has an empty interior.\\

One can also observe that $G$ acts  topologically freely iff  for
any finite set $\{ g_1 , ... , g_n \} \subset G \hspace{2mm}
(g_i \ne e)$ the set $[ \cup^n _{i = 1} X_{g_i} ]$ has an empty interior.\\

Just as in \cite{AnLeb}, 12.13 and 12.13' it can be shown that the
foregoing definition is equivalent to the next one: the action of
$G$ is said to be topologically free if for any finite set  $\{
g_1 , ... , g_k \} \subset G$ and a non empty open set $U\subset
X$ there exists a point $x\in U$ such that all the points $t_{g_i}
(x), \hspace{2mm}
i=1, ... , k$ are distinct.\\

Since $X$ is Hausdorff the latter definition is also equivalent to
the following: the action of $G$ is said to be topologically free
if for any finite set $\{ g_1 , ... , g_k \} \subset G$ and a non
empty open set $U\subset X$ there exists a non empty open set
$V\subset U$ such that
\begin{equation}
\label{e1.6'}
t_{g_i} (V) \cap t_{g_j} (V) =\emptyset \hspace{5mm}i,j \in \overline{1,k},\hspace{2mm}i\ne j.
\end{equation}
\end{I}
The next theorem is an analogue to Theorem 12.14 \cite{AnLeb}.
\begin{Tm}
\label{1.7}
If $G$ acts topologically freely then $B(A, T_g )$ possesses the property (*).
\end{Tm}
{\bf Proof:} Consider any element $b\in B(A, T_g )$ of the form
$$
b = \sum_{g \in F} a_g T_g
$$
where $F = \{ g_1 , ... , g_n \}.$ Fix any $\varepsilon > 0$ and
let $x_0 \in X$ be a point satisfying the condition
\begin{equation}
\label{e1.9}
\Vert a_e (x_0 ) \Vert > \Vert  a_e \Vert - \varepsilon
\end{equation}
(such a point $x_0$ exists since $A\cong C(X, L(E))$).\\
Consider a neighborhood $U$ of $x_0$ such that
\begin{equation}
\label{e1.10}
\Vert a_e (x_ ) \Vert >\Vert  a_e (x_0 ) \Vert -\varepsilon , \hspace{5mm} x\in U.
\end{equation}
In view of the topological freedom of the action of $G$ (see \ref{1.6}) there
exists a non empty open set
$V$ such that
\begin{equation}
\label{e1.11}
V\subset U \hspace{2mm}{\rm and}\hspace{2mm} t^{-1}_{g_i}(V) \cap  t^{-1}_{g_j}(V) =\emptyset ,
\hspace{5mm} i,j \in F, i\ne j.
\end{equation}
Let $c\in {\mathbb Z}(A) $ be a certain element having the form $c
= c(x)I$ where $c(\cdot )$ is a continuous function on $X$ with
the properties

(i) $0\le c(x) \le 1,$

(ii) $c(x^\prime ) = 1 $ for some $x^\prime \in V, $

(iii) $c(x)= 0,$ when $x \notin V.$\\
In view of (i) and (ii) and by (\ref{e1.9}) and  (\ref{e1.10}) and the choice of $x^\prime$
 we have
\begin{equation}
\label{e1.12}
\Vert c^2 a_e \Vert \ge \Vert a_e (x^\prime  ) \Vert > \Vert a_e \Vert - 2\varepsilon
\end{equation}
and
\begin{equation}
\label{e1.13}
\Vert cbc \Vert \le \Vert b \Vert .
\end{equation}
And (\ref{e1.11}) along with (iii) imply
\begin{equation}
\label{e1.14}
c{\hat T}_g (c) = 0, \hspace{5mm}g\ne e,\hspace{2mm} g\in F.
\end{equation}
Now by (\ref{e1.14}), (\ref{e1.13}) and (\ref{e1.12}) we obtain
$$
\Vert b \Vert \ge \Vert cbc \Vert = \Vert c (\sum a_g T_g )c\Vert = \Vert \sum c a_g {\hat T}_g (c) T_g \Vert =
 \Vert \sum c {\hat T}_g (c)a_g  T_g \Vert =
\Vert c^2 a_e \Vert  > \Vert a_e \Vert - 2\varepsilon .
$$
Which in view of the arbitrariness of $\varepsilon $ finishes the
proof.
\mbox{}\qed\mbox{}\\
\section{ Properties (*),  (**) and  the dual group action}
\setcounter{equation}{0}

We have already noted (Remark \ref{Rk1.4}) that in the $C^\ast -$algebra case the
property (**) (\ref{e1.6}) follows from the property (*). In fact the main reason for this
is that\\

 {\em the presence of the property (*) implies}  (\cite{AnLeb}, Theorem 12.8)
$$
C^\ast (A, T_g ) \cong A \times_{\hat T} G
$$
where by $ A \times_{\hat T} G$ we denote the cross-product of the algebra $A$ by the group
$\{ {\hat T}_g \}_{g \in G}$ of its automorphisms (here $G$ is considered as a discrete group).\\
Since in a Banach space case we do not have anything like the
isomorphism mentioned above we have to check the property (**)
even when $B(A, T_g )$ possesses the property (*). In a general
situation (that is for an {\em arbitrary} discrete group of
isometries $\{  T_g \}_{g \in G}$ with $T_g A T^{-1}_g = A$) the
verification of the property (**) may be very sophisticated.
Theorem \ref{1.13} presented below shows that in the case of a
{\em locally compact commutative} group $G$ and under a special
assumption (which as it will be seen later is in fact rather
common) the algebra
$B(A, T_g )$ possesses the properties (*) and (**) simultaneously.\\

Before proceeding to this theorem we would like to recall the
following {\em Kronecker approximation theorem} linking continuous
and discontinuous characters of locally compact commutative
groups.
\begin{Tm}
\label{1.12} {\rm (\cite{HewRoss}, (26.15)).} Let $G$ be a locally
compact commutative group, $\hat G$ be its dual group and $\Gamma$
be any subgroup of $G.$ Further let $f$ be any (not necessarily
continuous) character of $\Gamma$ and $g_1 , ... , g_n$ be some
elements of $\Gamma$. Then for any $\varepsilon > 0$ there exists
a {\em continuous} character $\chi \in {\hat G}$ such that
$$
\vert f(g_i ) - \chi (g_i )\vert < \varepsilon , \hspace{5mm}i ={\overline{1, n}}.
$$
\end{Tm}
\begin{Tm}
\label{1.13}
Let $G$ be a locally compact commutative group. If for any finite set $F\subset G$ and
any character $\chi \in {\hat G}$ there is satisfied the equality
\begin{equation}
\label{e1.36}
\Vert \sum_{g\in F} a_g T_g \Vert = \Vert \sum_{g\in F} a_g  \chi (g) T_g \Vert
\end{equation}
then the algebra $B(A, T_g )$ possesses the properties (*) and (**).
\end{Tm}
{\bf Proof:} In view of the preceding theorem if (\ref{e1.36})
holds for every $\chi \in {\hat G}$ it also holds for every not
necessarily continuous character of $G$ (that is the character of
$G$ considered as a discrete group). Bearing this in mind we {\em
henceforth assume  $G$ to be discrete and denote by $\hat G$
its dual as for the discrete group.}\\
Let $b \in B(A, T_g )$ be an element of the form
$$
b = \sum_{g\in G} a_g T_g  , \hspace{5mm} \vert F \vert < \infty .
$$
By $b(\chi ), \chi \in {\hat G}$ we denote the element of $B(A, T_g )$ given by
$$
b (\chi )= \sum_{g\in G} a_g \chi (g) T_g .
$$
(\ref{e1.36}) means that
\begin{equation}
\label{e1.37}
\Vert b \Vert = \Vert b(  \chi ) \Vert , \hspace{5mm} \chi \in {\hat G}.
\end{equation}
Let us prove first that this condition implies the property (*).\\

Recall that for any $c\in L(F)$ we have
\begin{equation}
\label{e1.38}
\Vert c \Vert = \sup_{\xi \in F, \Vert \xi \Vert = 1; \eta\in F^\ast  \Vert \eta \Vert = 1}
\vert <c\xi ,\eta > \vert
\end{equation}
here $F^\ast$ is the dual space of $F$ and $<\xi , \eta > = \eta (\xi )$ ($F$ is
the space where the
operators from $B(A, T_g )$ act).\\
Fix any $\varepsilon > 0$ and let $\xi _0 \in F$ and $\eta _0 \in F^\ast $, $\Vert \xi _0 \Vert =
\Vert \eta _0 \Vert = 1$ be such that
\begin{equation}
\label{e1.39}
\Vert a_e \Vert \ge \vert < a_e \xi_0 , \eta_0 > \vert > \Vert a_e \Vert - \varepsilon
\end{equation}
here $e$ is the identity of $G.$\\
Now (\ref{e1.37}) and (\ref{e1.38}) imply
$$
\Vert b \Vert  = \Vert  b (\chi ) \Vert  \ge  \vert < b (\chi ) \xi_0 , \eta_0 > \vert =
\vert < \sum_F a_g \chi (g) T_g \xi_0 , \eta_0 > \vert =
$$
\begin{equation}
\label{e1.40}
\vert  \sum_F <a_g T_g \xi_0 , \eta_0 >  \chi (g) \vert =
\vert \sum_F \alpha _g  \chi (g) \vert
\end{equation}
where $\alpha_g = < a_g T_g \xi_0 , \eta_0 >.$\\
Consider the function $f(\chi ) \in C({\hat G})$ given by
$$
f(\chi ) = \sum_F \alpha_g g(\chi )
$$
where $g(\chi ) = \chi (g)$ is the character on $\hat G$.\\
Since $G$ is discrete it follows that $\hat G$ is compact (\cite{HewRoss}, (23.17)) and
the Fourier coefficient $\alpha_e$ of $f$ can be calculated (recalling \cite{HewRoss}, (23.19))
as
\begin{equation}
\label{e1.41}
\alpha_e = \int_{\hat G} f(\chi ) d\mu
\end{equation}
($\mu$ is the normalized Haar measure on $\hat G$).\\
Which implies in turn (by the definition of $f$ and (\ref{e1.40}))
$$
\vert \alpha_e \vert \le \int_{\hat G} \vert f(\chi ) \vert d\mu \le
\max_{\hat G} \vert f(\chi ) \vert  \le \Vert b \Vert .
$$
And along with (\ref{e1.39}) this means that
$$
\Vert b \Vert \ge \Vert a_e \Vert - \varepsilon
$$
thus (in view of the arbitrariness of $\varepsilon $) proving the property (*).\\

Now we pass to the verification of the property (**).\\

Let $b$ be any element of $B(A, T_g ).$\\
If (\ref{e1.37}) holds then as we have already checked $B(A, T_g )$
possesses the property (*) and thus the mapping $N_{g_0 }$ given by (\ref{e1.5})
is correctly defined.\\
We must show that if $N_g (b) = 0$ for every $g \in G$ then $b =0$.\\
Evidently it is enough to prove that for any fixed $\chi \in F$ and $\eta \in F^\ast$
with $\Vert \xi \Vert = \Vert \eta \Vert = 1$ we have
\begin{equation}
\label{e1.42}
< b\xi , \eta > = 0.
\end{equation}
Let $\{ b_n \}$ be a sequence of elements of $B(A, T_g )$
tending to $b$ and each having the form
$$
b_n = \sum_{F_n} a^{(n)}_g T_g , \hspace{5mm} \vert F_n \vert < \infty .
$$
Consider the elements $b_n (\chi ), \chi \in {\hat G}$ given by
$$
b_n (\chi )= \sum_{F_n} a^{(n)}_g \chi (g) T_g
$$
and the sequence $\{  f_n (\chi ) \}$ of (continuous on $\hat G$) functions defined by
\begin{equation}
\label{e1.43}
 f_n (\chi ) = < b_n (\chi )\xi , \eta > = \sum_{F_n} < a^{(n)}_g T_g \xi , \eta > \chi (g) =
 \sum_{F_n} \alpha^{(n)}_g \chi (g)
\end{equation}
where $\alpha^{(n)}_g  = < a^{(n)}_g T_g \xi , \eta >$.\\
Since $b_n \to_{n\to\infty} b$ it follows (in view of (\ref{e1.37})) that
\begin{equation}
\label{e1.44}
 \Vert b_{n_1} (\chi ) - b_{n_2} (\chi )\Vert = \Vert b_{n_1} - b_{n_2} \Vert
 \to_{n_1 , n_2 \to \infty} 0
\end{equation}
which implies that for every fixed $\chi_0 \in {\hat G}$ the sequence $\{  b_n (\chi_0 )  \}$
tends to a certain element $ b (\chi_0 ) \in B(A, T_g )$.\\

Let $f(\chi )$ be the function given by
$$
f(\chi ) = < b(\chi )\xi , \eta  >.
$$
Applying (\ref{e1.44}) we obtain
$$
\vert   f_{n_1} (\chi ) - f_{n_2} (\chi ) \vert =
\vert < [ b_{n_1} (\chi ) - b_{n_2} (\chi ) ] \xi , \eta > \vert \le
\Vert b_{n_1} - b_{n_2} \Vert \to_{n_1 , n_2 \to \infty} 0.
$$
Which means that the sequence $\{ f_n  \}$ of (continuous) functions
tends uniformly (on $\hat G$)
to $f$. Thus $f$ is continuous and (since $\mu ({\hat G}) = 1$)
it follows that $f \in L^2 ({\hat G})$.\\
Let
$$
f(\chi ) = \sum_G \alpha_g g(\chi )
$$
where  the righthand part is the Fourier series of $f.$\\
Since $f_n \to f $ (in $ L^2 ({\hat G})$) it follows that
\begin{equation}
\label{e1.45}
 \alpha^{(n)}_g \to \alpha _g \hspace{2mm} {\mbox {\rm for every}}\hspace{2mm} g\in G
\end{equation}
where $\alpha^{(n)}_g $ are those defined by (\ref{e1.43}).\\
Now note that the property (*) implies
\begin{equation}
\label{e1.46}
 \Vert a^{(n)}_g  \Vert \to \Vert N _g (b) \Vert \hspace{2mm} {\mbox {\rm for every}}\hspace{2mm}
 g\in G.
\end{equation}
And also observe  that
$$
\vert \alpha^{(n)}_g  \vert= \vert < a^{(n)}_g  T_g \xi , \eta >  \vert \le \Vert a^{(n)}_g  \Vert
$$
which together with (\ref{e1.45}) and (\ref{e1.46}) means that
\begin{equation}
\label{e1.47}
\alpha_g = 0  \hspace{2mm} {\mbox {\rm for every}}\hspace{2mm} g\in G.
\end{equation}
And (\ref{e1.47}) and the continuity of $f$ implies
$$
f(\chi ) = 0 \hspace{2mm} {\mbox {\rm for every}}\hspace{2mm} \chi\in{\hat G}.
$$
In particular
$$
f(1) = < b\xi , \eta > = 0.
$$
Thus (\ref{e1.42}) is true and the proof is finished.
\mbox{}\qed\mbox{}\\
\begin{Rk}
\label{1.14}
\em
(1) The theorem just proved asserts that (in the case of a commutative group $G$)
the property (\ref{e1.36}) implies the properties (*) and (**). On the other hand in the $C^\ast -$algebra
case the property (*) implies the property  (\ref{e1.36}). Indeed by \cite{AnLeb}, Theorem 12.8
if $C^\ast (A, T_g )$ possesses the property (*) then
$$
C^\ast (A, T_g ) \cong A\times_{\hat T} G
$$
and the isomorphism is established by the mapping
\begin{equation}
\label{e1.48}
\sum a_g T_g \longleftrightarrow \sum a_g \otimes g
\end{equation}
If for any fixed character $\chi \in {\hat G}$ we define the representation
$$
\pi_\chi : A\times_{\hat T} G \to A\times_{\hat T} G
$$
by the formulae
\begin{equation}
\label{e1.49}
\pi_\chi (  \sum a_g \otimes g ) =  \sum \chi (g)a_g \otimes g
\end{equation}
then
\begin{equation}
\label{e1.50}
\Vert \pi_\chi (  \sum a_g \otimes g ) \Vert  = \Vert \sum a_g \otimes g \Vert
\end{equation}
(which follows again from  \cite{AnLeb}, Theorem 12.8 and the observation that
$\pi_\chi (A\times_{\hat T} G  ) $  possesses the property (*) since $A\times_{\hat T} G $
possesses it (see  \cite{AnLeb}, 12.3)).\\
Now (\ref{e1.50}) along with (\ref{e1.49}) and    (\ref{e1.48}) means that
$$
\Vert \sum a_g \chi (g) T_g \Vert = \Vert \sum a_g  T_g \Vert
$$
and consequently (\ref{e1.36}) is true.\\
Thus in the $C^\ast -$algebra case (and when $G$ is commutative and locally compact)
 the properties (*)
and  (\ref{e1.36}) are equivalent.\\
So in the general Banach algebra case one can regard the condition
(\ref{e1.36}) as a certain substitute for the property (*) (under
which the properties (*) and (**) are realized simultaneously as
in the $C^\ast -$case).\\

(2) It is worth mentioning that if $G$ is not commutative then one  has to consider
the conditions different from   (\ref{e1.36})  (recall for example  \cite{AnLeb}, Lemma 22.6).\\

(3) We would also like to note that while in Theorem \ref{1.7}  we assumed the algebra
$A$  to have a special form (described in
\ref{1.5}) the algebra $A$ considered in Theorem  \ref{1.13} is an {\em arbitrary} Banach algebra.
\end{Rk}
Now we present an example of the algebra $B(A, T_g )$ for which the verification
of the property (\ref{e1.36}) could be carried out in an easy way.
\begin{Ee}
\label{1.15}
\em
Let $X = G$ be a locally compact commutative group, $A$ be a certain algebra of
 operators isomorphic to
$C(X, L(E))$ and $\{ T_g \}_{g\in G}$ be a group of isometries satisfying (\ref{e1.1})
and such that
$$
(T_g a T^{-1}_g )(x) = a (x\cdot g), \hspace{5mm} a\in A.
$$
In this case \\

{\em the algebra $B(A, T_g )$ generated by $A$ and $\{ T_g \}_{g\in G}$ possesses
the property (\ref{e1.36}).  }\\

Indeed. For every $\chi \in {\hat G}$ consider the operator ${\bar \chi} \in {\mathbb Z}(A)$
corresponding to the
operator valued function $\chi (x)I$ ($I$ is the identity operator in $E$).
Obviously ${\bar \chi}$ is an invertible isometry and for any operator
$$
b = \sum_{g\in F} a_g T_g , \hspace{5mm} \vert F \vert < \infty
$$
we have
$$
\Vert b \Vert = \Vert {\bar \chi}^{-1} b {\bar \chi} \Vert = \Vert \sum a_g \chi (g) T_g \Vert .
$$
Thus  (\ref{e1.36}) is true.
\mbox{}\qed\mbox{}\\
\end{Ee}
\begin{Rk}
\label{1.16}
\em
In fact here we have not used the commutativity of $G$ and the analogous
statement can be proved for any topological group.
\end{Rk}
\begin{Ee}
\label{1.17}
\em
{ \bf Regular representation of an algebra and a group of automorphisms. }
Let us also consider one more example where the properties (*) and (**) can be checked easily
--- {\em the regular representation of an algebra $A$ and a group of automorphisms
$\{ {\hat T}_g \}_{g\in G}$. }\\
Namely let $A\subset L(D)$ be a certain Banach algebra and $\{ {\hat T}_g \}_{g\in G}$
be a certain
group of its automorphisms ($G$ is an {\em arbitrary} group that is
{\em not necessarily commutative}).\\
Denote by $H$ any of the spaces $l^p (G, D), 1\le p \le \infty$ or $l_0 (G, D)$
(here $l_0 (G, D)$ is the space of vector valued functions on $G$ having values
in $D$ and tending to zero at infinity (with the sup-norm)).\\
Set the operators $V_{g_0} : H \to H$ by the formula
\begin{equation}
\label{e1.51}
(V_{g_0} \xi )(g) = \xi (gg_0 ), \hspace{5mm} g, g_0 \in G
\end{equation}
and consider the algebra ${\bar A} \subset L(H)$ isomorphic (as a Banach algebra)
to $A$ and given by
\begin{equation}
\label{e1.52}
({\bar a}\xi )(g) = {\hat T}_g (a) \xi (g), \hspace{5mm} a \in A.
\end{equation}
Routine computation shows that with this notation we have
$$
V_g {\bar a} V^{-1}_g = {\overline{{\hat T}_g (a)}}
$$
which in view of the isomorphism between $A$ and $\bar A$ means that the operators $V_g , g\in G$
given by (\ref{e1.51}) generate the automorphisms ${\hat T}_g$ of $\bar A$.\\

The algebra $B({\bar A}, V_g )\subset L(H)$ is called the {\em (right) regular representation}
corresponding to the algebra $A$ and the group of automorphisms $\{ {\hat T}_g \}_{g\in G}$
(in fact we have the series of representations depending on the type of the space
$H$ chosen).
\end{Ee}
We observe that \\

{\em the algebra $B({\bar A}, V_g )$ possesses the properties (*),  (**) and  (\ref{e1.36})
(for every $H$ considered)}.\\

Let us check  the property (*) first.\\
Take any operator $b$ of the form
\begin{equation}
\label{e1.53}
b = \sum_{g \in F} {\bar a} V_g , \hspace{5mm} \vert F \vert < \infty .
\end{equation}
Fix  $\varepsilon > 0$ and choose  $\eta \in D , \hspace{2mm} \Vert \eta \Vert = 1$
such that
$$\Vert a_e \eta \Vert \ge \Vert a_e \Vert - \varepsilon .$$
 Consider the vector $\xi \in H$ defined by
$$
\xi (g) = \delta (e,g) \eta
$$
where $ \delta (e,g)$ is the Kronecker symbol. From the definition of $\xi$ and the
explicit form of $b$
we have
$$
\Vert b \xi \Vert = \Vert  \sum_{g \in F} {\bar a} V_g (\xi ) \Vert \ge
\Vert {\bar a}_e \xi \Vert =
\Vert a_e \eta \Vert \ge \Vert a_e \Vert - \varepsilon =  \Vert {\bar a}_e \Vert - \varepsilon
$$
which (in view of the arbitrariness of $\varepsilon$) proves the property (*).\\

Now let us check the property (**).\\
We start with the case when $H$ is equal to $l^p (G, D), 1\le p < \infty$ or $l_0 (G, D)$.\\
In this situation any operator $b \in L(H)$ is defined by its matrix
$$
[b]_{h,g} = {\hat T}_h [N (b V_{g^{-1}h})], \hspace{5mm} h,g \in G
$$
(here $N : B({\bar A}, V_g ) \to A\cong {\bar A}$ is defined by  (\ref{e1.5}) and
the isomorphism between
$A$ and ${\bar A}$).\\
Since an operator is zero iff its matrix is zero the property (**) for
 the situation considered is proved.\\

Now it remains to consider the case $H = l^\infty (G, D)$.\\
Here observe first that for any operator $b$ of the form (\ref{e1.53}) we have
\begin{equation}
\label{e1.54}
\Vert b \Vert_{l^\infty (G, D)} =  \Vert b \Vert_{l_0  (G, D)}.
\end{equation}
Indeed since $l_0  (G, D) \subset l^\infty (G, D)$  it follows that
$$
\Vert b \Vert_{l^\infty (G, D)} \ge   \Vert b \Vert_{l_0  (G, D)}.
$$
To obtain the opposite inequality
 consider for a fixed $\varepsilon > 0$ any vector $\xi \in l^\infty (G, D)$ with the properties
\begin{equation}
\label{e1.55}
\Vert \xi \Vert =  1
\end{equation}
and
\begin{equation}
\label{e1.56}
\Vert b \xi \Vert_{l^\infty (G, D)} > \lambda - \varepsilon
\end{equation}
where $\lambda = \Vert b \Vert_{l^\infty (G, D)}$.\\
The property (\ref{e1.56}) means that there exists a point $g_0 \in G$ such that
\begin{equation}
\label{e1.57}
\Vert (b \xi ) (g_0 ) \Vert > \lambda - \varepsilon .
\end{equation}
Set the vector $\eta \in l_0 (G, D)$ by the equations
$$
\eta (g) =
\left \{ \begin{array}{ll}
\xi (g), \hspace{2mm} g\in M = [{\cup_{g\in F}} (gg_0 )]\\
0,\hspace{2mm} g\notin M.
\end{array} \right.
$$
For this $\eta $ in view of (\ref{e1.55}) and (\ref{e1.57}) we have
$$
\Vert \eta \Vert_{l_0 (G, D)} \le 1
$$
and
$$
\Vert b \eta \Vert_{l_0 (G, D)}>   \lambda - \varepsilon
$$
which (as $\varepsilon $ is arbitrary) implies
$$
\Vert b \eta \Vert_{l_0 (G, D)}\ge   \lambda =  \Vert b  \Vert_{l^\infty (G, D)}
$$
thus proving (\ref{e1.54}).\\
Since the equality  (\ref{e1.54}) holds for any $b$ having the form  (\ref{e1.53})
it also holds for any
$b \in B({\bar A}, V_g )$. Thus as $B({\bar A}, V_g )$ possesses the property (**) as acting
in $l_0 (G, D)$ it possesses this property as acting in $l^\infty (G, D) $ as well.\\

Finally we have to check the equality (\ref{e1.36}).\\
Take any character $\chi \in {\hat G}$ and let ${\bar \chi} \in L(H)$ be the isometry
 defined by the operator
valued function $\{   \chi (g)I \}_{g \in G}$ ($I$ is the identity operator in $D$).
Obviously ${\bar \chi}$
commutes with every element ${\bar a} \in {\bar A}$ and for any operator $b$ of
 the form (\ref{e1.53})
we have $$
\Vert b \Vert = \Vert   {\bar \chi}^{-1} b  {\bar \chi}  \Vert  = \Vert \sum_{g\in F}
{\bar a}_g
 \chi (g) V_g \Vert
$$
which is the desired result.
\begin{Rk}
\label{1.18}
\em
(1) We would like to emphasize once more that since here $G$ is not assumed
 to be commutative the equality
(\ref{e1.36}) does not imply the properties (*) and (**).\\
(2) Example \ref{1.17} presents a very important for our subsequent purposes class
of algebras (the regular representations) possessing the properties (*) and (**).
The essence of the argument used in this example shows that the main point why
the properties (*) and (**) are satisfied here is that the regular representations
consists of operators acting in {\em discrete spaces} ($l^p (G, D), \hspace{2mm} l_0 (G, D)$)
 and consequently such operators can be presented by their matrices.\\
Namely we have in mind the following. Let $H$ be any of the spaces $l^p (G, D),
 \hspace{2mm} 1\le p < \infty$ or
$l_0 (G, D)$ where $D$ is a certain Banach space. Let ${\sf A} \subset L(H)$
be the algebra of operators of the form
$$
(a \xi )(g) = a (g) \xi (g) , \hspace{5mm}\xi \in H, g \in G, a \in {\sf A}
$$
where $a (g) \in L(D)$ and $\sup_G  \Vert a(g) \Vert < \infty$
(thus $ {\sf A} \cong l^\infty  (G, L(D)) $).\\
Let the operators ${\sf V}_{g_0} : H \to H$ be given by the formula
$$
({\sf V}_{g_0} \xi) (g) = \xi (g g_0 ), \hspace{5mm} g , g_0 \in G.
$$
Consider the algebra $B ({\sf A}, {\sf V}_g ) \subset L(H)$.
The obvious modification of the argument used for the consideration of $B({\bar A}, V_g )$
in Example \ref{1.17} shows that\\

{\em $B ({\sf A}, {\sf V}_g ) $ possesses the properties (*), (**) and (\ref{e1.36})
(for every $H$ considered). }\\

In fact the algebra $B ({\sf A}, {\sf V}_g ) $ in connection with
the properties examined 'dominates'
$B({\bar A}, V_g )$ since for any fixed $H$ (mentioned in Example  \ref{1.17}) we have
$$
{\bar A} \subset {\sf A}
$$
and
$$
V_g = {\sf V}_g
$$
which means that in this case $B({\bar A}, V_g )$ is a subalgebra of $B ({\sf A}, {\sf V}_g ) $.
\end{Rk}
Now  we shall present a series of examples of Banach algebras generated by
'weighted composition operators' acting in various spaces. We shall study the structure of these
algebras and find out in particular
that
in a number of
 natural situations the condition on $G$ to act
topologically freely implies the satisfaction of the property (**) (condition (\ref{e1.36})).\\

In all the cases we shall specify not only the form of the algebra $A\subset L(D)$ but
also the space $D$ and the explicit form of the operators $T_g , g\in G$.
This is done in view of the fact that the norm
$
\Vert \sum a_g T_g \Vert
$
 is calculated (estimated)  in different spaces
 in different ways even when the operators $T_g$ generate the same automorphisms of $A$.\\
\section{ Example 1. Operators in $C(X, E)$}
\setcounter{equation}{0}

\begin{I}
\label{2.1}
\em
 Let $D =C(X, E)$ be the Banach space of all continuous functions on a completely regular space
$X$ taking values in a Banach space $E$ with the sup-norm.\\
Consider the algebra $A = C(X, L(E)) \subset L(D)$ of operators of multiplication by
continuous operator valued functions with the sup-norm, that is for any $a\in A$ and $f\in D$
\begin{equation}
\label{e2.1}
(af)(x) = a(x) f(x) .
\end{equation}

Henceforth when considering the operators of the form (\ref{e2.1})
(where $f = \{ f(x) \}$ is a certain vector valued function and $a
= \{ a(x) \}$ is a certain operator valued function  ($a(x)\in
L(F), \hspace{2mm}F\ni f(x)$)) we shall call them the
{\em multiplication} (or  {\em diagonal}) operators.\\

Let $\{  t_g \}_{g \in G}$ be a group of homeomorphisms of $X$.
By $T_g$ we denote the isometry of $D$ defined by the formula
\begin{equation}
\label{e2.2}
(T_g f)(x) =  f( t^{-1}_g x) .
\end{equation}
Evidently $T_g$ satisfies (\ref{e1.1}) and the homeomorphisms $t_g$ mentioned above
coincide with those
given by (\ref{e1.7}).
\end{I}
For the algebra $B(A, T_g ) \subset L(D)$ generated by these $A$ and $T_g$ we can
calculate even more than the properties (*) and (**) and the equality (\ref{e1.36}).
Namely let $F\subset G$ be any fixed finite set, we denote by $B_F (E)$ and $S_F (E)$
respectively the sets
\begin{equation}
\label{e2.3}
B_F (E) = \{  \{  f_g  \}_{g\in F} : f_g \in E, \hspace{2mm}\Vert f_g \Vert \le 1,
 \hspace{2mm}g\in F  \} ,
\end{equation}
\begin{equation}
\label{e2.4}
S_F (E) = \{  \{  f_g  \}_{g\in F} : f_g \in E, \hspace{2mm}\Vert f_g \Vert = 1,
\hspace{2mm}g\in F  \} .
\end{equation}
We start the consideration of the algebra $B(A, T_g )$ with the
following observation.
\begin{Tm}
\label{2.2}
If $G$ acts topologically freely then
$$
\Vert \sum_{g\in F} a_g T_g \Vert = {\sup_X} \sup_{  \{  f_g  \}_{g\in F} \in S_F (E) }
\Vert \sum_{g\in F} a_g  (x) f_g  \Vert =
$$
\begin{equation}
\label{e2.5}
 {\sup_X} \sup_{  \{  f_g  \}_{g\in F} \in B_F (E) }
\Vert \sum_{g\in F} a_g  (x) f_g  \Vert  .
\end{equation}
\end{Tm}
{\bf Proof:} To shorten the notation we omit the space symbol $(E)$ in $B_F (E) $ and $S_F (E)$.\\
We start with the establishing of the equality
\begin{equation}
\label{e2.6}
\Vert \sum_{g\in F} a_g T_g \Vert =
 {\sup_X} \sup_{  \{  f_g  \}_{g\in F} \in B_F (E) }
\Vert \sum_{g\in F} a_g  (x) f_g  \Vert  .
\end{equation}
Observe that for any $f\in C(X, E)$ with $\Vert f \Vert =1$ we have
$$
\Vert \sum_{g\in F} a_g T_g f \Vert =
 {\sup_X} \Vert  \sum_{g\in F} a_g (x) f(t^{-1}_g (x)) \Vert =
 {\sup_X} \Vert  \sum_{g\in F} a_g (x) f_g (x) \Vert
$$
where $\Vert  f_g (x) \Vert   = \Vert  f(t^{-1}_g (x)) \Vert \le 1$.\\
And it follows that
$$
\Vert \sum_{g\in F} a_g T_g  \Vert \le  {\sup_X} \sup_{  \{  f_g  \} \in B_F  }
 \Vert  \sum_{g\in F} a_g (x) f_g \Vert .
$$
to establish the opposite inequality fix any $x_0 \in X$ and any collection
$ \{  f_g  \} \in B_F $ and let
\begin{equation}
\label{e2.7}
\Vert  \sum_{g\in F} a_g (x_0 ) f_g \Vert =\lambda .
\end{equation}
Since $a_g \in C(X, L(E))$ for any $\varepsilon > 0$ there exists
a neighborhood $U$ of $x_0$ such that
\begin{equation}
\label{e2.8}
\Vert  \sum_{g\in F} a_g (x ) f_g \Vert  > \lambda -\varepsilon
\hspace{2mm} {\mbox {\rm for every}}\hspace{2mm} x\in U.
\end{equation}
as $G$ acts topologically freely one can find a point $x^\prime
\in U$ and its neighborhood $V\in U$ such that
\begin{equation}
\label{e2.9}
V \cap t^{-1}_g (V) = \emptyset ,\hspace{5mm} g\in F, g\neq e
\end{equation}
Take some  functions $\phi _g \in C(X), \hspace{2mm} g\in F$ with the properties
\begin{equation}
\label{e2.10}
\left \{ \begin{array}{lll}
0 \le \phi_g (x) \le 1\\
\phi_g (t^{-1}_g (x^\prime )) = 1, \hspace{5mm} g\in F\\
\phi_g (X\setminus t^{-1}_g (V )) =), \hspace{5mm} g\in F
\end{array} \right.
\end{equation}
and set
\begin{equation}
\label{e2.11}
f = \sum_{g\in F} \phi_g (x)f_g \in C(X, E) .
\end{equation}
Now by the choice of $\phi_g$ and since
$\{ f_g  \} \in B_F$ we have
\begin{equation}
\label{e2.12}
\Vert f \Vert \le 1 .
\end{equation}
And on the other hand
\begin{equation}
\label{e2.13}
\Vert (\sum a_g T_g f) (x^\prime ) \Vert =
 \Vert  \sum a_g (x^\prime ) f_g \Vert   \ge \lambda -\varepsilon
\end{equation}
where the latter inequality follows from (\ref{e2.8}).\\
Thus by the definition of $\lambda $ (see (\ref{e2.7})) and in
view of the arbitrariness of $\varepsilon$ and (\ref{e2.12}) the
inequality (\ref{e2.13}) means that
$$
\Vert  \sum a_g T_g  \Vert \ge  {\sup_X} \sup_{  \{  f_g  \} \in B_F  }
 \Vert  \sum_{g\in F} a_g (x) f_g \Vert .
$$
which finishes the proof of (\ref{e2.6}).\\

Observe now that (\ref{e2.5}) follows from  (\ref{e2.6}).\\
Indeed.\\
 Consider the space
$$
{\tilde E} = E_1 \times E_2 \times ... \times E_{\vert F \vert} ,\hspace{5mm} E_i = E
$$
with the norm
$$
\Vert {\tilde e}\Vert_{\tilde E}= \Vert (e_1 , ... , e_{\vert F \vert} ) \Vert_{\tilde E} =
\max_{i {\overline{1, \vert F \vert}}} \Vert e_i \Vert_E
$$
and the operator ${\tilde b}(x) \in L({\tilde E}, E)$ given by
$$
{\tilde b}(x) ( (e_1 , ... , e_{\vert F \vert} ) ) =
a_{g_1 }(x)e_1 + ... +  a_{g_{ \vert F \vert} }(x)e_{ \vert F \vert}
$$
($\{ g_1 , ... , g_{\vert F \vert}   \} = F, \hspace{2mm} x\in X$).\\
Obviously $B_F (E)$ can be identified with the unit ball of ${\tilde E}$ and $B_F (E)$ is
the convex hull of $S_F (E)$.\\
Since the function
$$
\phi ({\tilde e}) = \Vert {\tilde b}(x)({\tilde e})\Vert_{\tilde E}
$$
is a convex function on $ {\tilde E}$ the  interrelation between  $B_F (E)$
and $S_F (E)$  mentioned above implies
$$
 \Vert {\tilde b}(x) \Vert = \sup_{{\tilde e} \in
  B_F (E) } \Vert {\tilde b}(x) ({\tilde e} )\Vert =
\sup_{{\tilde e} \in  S_F (E) } \Vert {\tilde b}(x) ({\tilde e} )\Vert  .
$$
Which along with the explicit form of ${\tilde b}(x)  $ means that
(\ref{e2.6}) and (\ref{e2.5}) are equivalent.
\mbox{}\qed\mbox{}\\
\begin{Rk}
\label{2.3}
\em
(1) If $E ={\mathbb C}$ (that is $D = C(X)$ and $A = C(X)$) then (\ref{e2.5}) implies\\

{\em if $G$ acts topologically freely then }
\begin{equation}
\label{e2.21}
\Vert  \sum_F a_g T_g  \Vert =  {\sup_X}
  \sum_F  \vert a_g (x)  \vert  .
\end{equation}
Indeed on the one hand
$$
{\sup_X}  \sup_{  \{  f_g  \}_{g\in F} \in S_F ({\mathbb C} ) } \vert  a_g (x) f_g  \vert  \le
{\sup_X}
  \sum  \vert a_g (x)  \vert
$$
and to obtain the opposite inequality just set for every $x\in X$ and $g\in F$
$$
f_g (x) =
\left \{ \begin{array}{ll}
[\arg a_g (x)]^{-1}, \hspace{2mm} {\mbox {\rm if}}  \hspace{2mm} a_g (x) \neq 0 \\
1,\hspace{2mm} {\mbox {\rm if}}  \hspace{2mm} a_g (x) = 0 .\\
\end{array} \right.
$$
(2) The equality (\ref{e2.5}) also shows that\\

{\em if $G$ acts topologically freely then}
$$
\Vert  \sum_F a_g T_g  \Vert = \Vert {\tilde b}_F \Vert
$$
{\em where}
$$
 {\tilde b}_F : D_1 \times D_2 \times ... \times D_{\vert F \vert} \to D, \hspace{5mm} D_i = D
$$
{\em is given by}
\begin{equation}\label{e2.22}
{\tilde b}_F (\xi_1 , ... , \xi_{\vert F \vert}) = a_{g_1}\xi_1 + ...
+ a{g_{\vert F \vert} }\xi_{\vert F \vert}
\end{equation}
($\{ g_1 , ... , g_{\vert F \vert}  \} = F$).
\end{Rk}
Now we shall continue the investigation of the algebra $B(A, T_g )$
and find out that in fact the equality
 (\ref{e2.5}) also leads to the next
\begin{La}
\label{2.4}
Let $B(A, T_g )$ be the algebra described  in 2.1 and $B({\bar A}, V_g )$ be the
corresponding regular representation for $A$ and $\{ {\hat T}_g \}_{g \in G}$
in the space $H = l_0 (G, C(X, E))$ (or  $ l^{\infty} (G, C(X, E))$). If $G$ acts
topologically freely then the mappings
$$
a \to {\bar a}, \hspace{5mm}a\in A
$$
$$
T_g \to V_g, \hspace{5mm}g\in G
$$
establish the isomorphism
$$
B(A, T_g ) \cong B({\bar A}, V_g ).
$$
\end{La}
{\bf Proof:} The simple calculation shows that the norm of the element
$$
{\bar b} = \sum_{g \in F} {\bar a}_g V_g \in B({\bar A}, V_g )
$$
is equal to
$$
{\sup_X} \sup_{  \{  f_g  \} \in B_F (E) }
 \Vert  \sum_{g\in F} a_g (x) f_g \Vert
$$
which in view of (\ref{e2.5}) finishes the proof.
\mbox{}\qed\mbox{}\\
\begin{Cy}
\label{2.5}
If $B(A, T_g )$ is that considered in 2.1 and $G$ acts topologically freely then
$B(A, T_g )$ possesses the properties (*) and (**) and (\ref{e1.36}).
\end{Cy}
{\bf Proof:} Follows from Lemma \ref{2.4} along with Example \ref{1.17}.
\mbox{}\qed\mbox{}\\
\section{ Example 2.  Operators in $L^{p}_\mu (X, E),$ ${\mbox 1< p < \infty}$}
\setcounter{equation}{0}

\begin{I}
\label{2.6}
\em
 Let $X$ be a completely regular space and $\mu$ be a certain Borel $\sigma -$finite
measure on $X$ with the support equal to the whole of $X$. Consider the Banach space
$D = L^{p}_\mu (X, E)$ where $p$ is a certain number $1<p< \infty$ (the cases $p = 1$ and
$p = \infty$ will be considered separately in Examples 3 and 4).\\

Let $A= C(X, L(E))$ be the algebra of the operators defined by (\ref{e2.1}) and let
 $\{ {\tilde t}_g \}_{g \in G}$ be a group of homeomorphisms of $X$ preserving the
equivalence class of $\mu$ (that is $t_g (\mu )$ is absolutely continuous with respect to $\mu$).
We denote by $T_g$ the isometry of $D$ given by
\begin{equation}
\label{e2.23}
(T_g f)(x) = \left[\frac{t^{-1}_g (\mu )} {d\mu}\right]^{1\over p} f(t^{-1}_g (x ))
\end{equation}
where $\frac{t^{-1}_g (\mu )} {d\mu}$ is the Radon-Nikodim derivative of $t^{-1}_g (\mu )$
with respect to $\mu$.\\
Let $B(A, T_g )\subset L(D)$ be the algebra generated by $A$ and   $\{ {\hat T}_g \}_{g \in G}$.
\end{I}
Our investigation of $B(A, T_g )$ will be divided into a series of steps and statements.\\
First of all we need to introduce a certain representation of  $B(A, T_g )$
 related to the trajectories
of points in $X$.\\
Namely for every point $x\in X$ we define the representation
\begin{equation}
\label{e2.24}
\pi_x : B(A, T_g ) \to L(l^p (G, E))
\end{equation}
$$
\pi_x (b) =b_x , \hspace{5mm}b\in B(A, T_g )
$$
by the equations
\begin{equation}
\label{e2.25}
(\pi_x (a)\xi )_g = a(t^{-1}_g (x))\xi_g ,
\end{equation}
\begin{equation}
\label{e2.26}
(\pi_x (T_{g_0 })\xi )_g = \xi_{g g_0}
\end{equation}
where $\xi = (\xi_g  )_{g \in G} \in l^p (G, E)$ and $a \in A$.\\
With this notation our first observation is
\begin{La}
\label{2.7}
 Let $b \in B(A, T_g ) \subset L(L^{p}_\mu (X, E))$
be an operator of the form
\begin{equation}
\label{e2.27}
b = \sum_F a_g T_g , \hspace{5mm} \vert F \vert < \infty
\end{equation}
($B(A, T_g )$ is that described in \ref{2.6}).\\
If $G$ acts topologically freely then
\begin{equation}
\label{e2.28}
\Vert b \Vert \ge \sup_x \Vert b_x \Vert
\end{equation}
where $b_x$ are given by (\ref{e2.24})-(\ref{e2.26}).\\
(Thus in particular $b \to b_x$ is indeed a representation of $B(A, T_g )$).
\end{La}
{\bf Proof:}  Fix some $x\in X$ and $\varepsilon > 0$ and let $\eta \in l^p (G, E)$
be a certain vector with the properties
\begin{equation}
\label{e2.29}
\Vert \eta \Vert = 1
\end{equation}
and
\begin{equation}
\label{e2.30}
\Vert b_x \eta \Vert \ge \Vert b_x  \Vert -\varepsilon
\end{equation}
where $b$ is given by (\ref{e2.27}). Without the loss of generality
 we can assume that only a finite
number of coordinates $\eta_g$ of $\eta = (\eta_g )$ are non zero that is there exists
a finite set $M\subset G$ satisfying
\begin{equation}
\label{e2.31}
\eta_g = 0 \hspace{2mm}{\mbox {\rm when}} \hspace{2mm} g \notin M.
\end{equation}
Let
$$
F_1= \bigcup_{g\in F} M\cdot g^{-1}.
$$
Clearly $\vert F_1 \vert < \infty$.\\
Since $G$ acts topologically freely it follows that there exists a non empty
open set $V\subset X$ such that $\mu (V) >0$ and
\begin{equation}
\label{e2.32}
t^{-1}_{g_1} (V) \cap t^{-1}_{g_2} (V) =\emptyset ,\hspace{5mm} g_1 , g_2 \in
[F_1 \cup M], \hspace{5mm}g_1 \neq g_2
\end{equation}
and for every $g_0 \in F$ and $g_1 \in F_1$ it is held
\begin{equation}
\label{e2.33}
\Vert a_{g_0}(t^{-1}_{g_1} (x)) -  a_{g_0}(y) \Vert < \varepsilon ,
\hspace{5mm} y \in  t^{-1}_{g_1} (V)
\end{equation}
(just choose $V$ in (\ref{e1.6'}) being contained in a
sufficiently small neighborhood $U$ of $x$ and observe that $\mu
(V) > 0$ as
${\rm supp} \mu = X$).\\
Let $\tilde T_g$ be the operator acting in $L^{p}_\mu (X) $
and defined by the formula (\ref{e2.23}). Set
${\bar \eta} \in L^{p}_\mu (X, E) $ by the formula
\begin{equation}
\label{e2.34}
{\bar \eta} =
\sum_{g\in M}\left[ (\mu (V))^{- {1\over p}}{\tilde T}_{g^{-1}}(\chi_V  )\eta_g \right] =
 \sum_{g\in M}{\bar \eta}_g
\end{equation}
where $\chi_V$ is the characteristic function of $V$.\\
The explicit form of $\eta$ along with (\ref{e2.32}), (\ref{e2.31}) and (\ref{e2.29})
implies
\begin{equation}
\label{e2.35}
  \Vert  {\bar \eta}_g \Vert  = \Vert   \eta_g \Vert  \hspace{5mm}{\rm and}  \hspace{5mm}
 \Vert  {\bar \eta} \Vert  = \Vert   \eta \Vert =1 .
\end{equation}
And in view of (\ref{e2.25}) and  (\ref{e2.26}) along with  (\ref{e2.31})  we obtain
 \begin{equation}
\label{e2.36}
(b_x \eta )_g = 0  \hspace{5mm}{\rm when}  \hspace{5mm}g\notin F_1
\end{equation}
and for every $g^{\prime} \in F_1$ (\ref{e2.25}) and  (\ref{e2.26}) give
\begin{equation}
\label{e2.37}
(b_x \eta )_{g^\prime} = \sum_{g\in F} a_g (t^{-1}_{g^\prime} (x))\eta_{{g^\prime}\cdot g}
\end{equation}
Considering now the vector $\bar \eta$ we have
\begin{equation}
\label{e2.38}
(b {\bar \eta} ) (y) = 0  \hspace{5mm}{\mbox {\rm for every}}  \hspace{5mm}
y\notin \cup_{g\in F_1}t^{-1}_g (V)
\end{equation}
and for every $g^{\prime} \in F_1$ and $ y \in t^{-1}_{g^\prime} (V)$
the explicit form of $\bar \eta$  (\ref{e2.34}) and $T_g$ (\ref{e2.23}) gives us
$$
(b {\bar \eta})(y) =  \sum_{g\in F} a_g (y) T_g {\bar \eta}_{{g^\prime}\cdot g} =
[\mu (V) ]^{-{1\over p}}
\sum_{g\in F} a_g (y) T_g ({\tilde T}_{(g^{\prime} g)^{-1}}(\chi_V  )\eta_{g^{\prime} g}  ) =
$$
$$
[\mu (V) ]^{-{1\over p}}
\sum_{g\in F} a_g (y) {\tilde T}_{{g^{\prime}}^{-1}}(\chi_V  )\eta_{g^{\prime} g}
$$
and it follows that
\begin{equation}
\label{e2.39}
\Vert   \chi_{t^{-1}_{g^\prime} (V)}\cdot b {\bar \eta} \Vert^p =\int_{t^{-1}_{g^\prime} (V)}
\Vert  \sum_{g\in F} a_g (y) \eta_{g^{\prime} g}   \Vert^p
\frac{\left[ {\tilde T}_{{g^{\prime}}^{-1}}(\chi_V  )  \right]^p}{\mu (V)} d\mu .
\end{equation}
Now recalling that ${\tilde T}_{{g^{\prime}}^{-1}}$ is an isometry in $L^{p}_\mu (X)$
and consequently
\begin{equation}
\label{e2.40}
\int_{t^{-1}_{g^\prime} (V)}\frac{\left[ {\tilde T}_{{g^{\prime}}^{-1}}(\chi_V  )  \right]^p}
{\mu (V)} d\mu
= 1
\end{equation}
we conclude (having in mind (\ref{e2.33})) that (\ref{e2.37}), (\ref{e2.39}) and (\ref{e2.40})
means that choosing $\varepsilon$ in (\ref{e2.33}) to be sufficiently small
we can obtain for any $\delta$
$$
\vert \Vert  (b_x \eta )_{g^\prime } \Vert^p - \Vert  \chi_{t^{-1}_{g^\prime} (V)} b{\bar \eta}
 \Vert^p    \vert
<\delta\hspace{5mm} {\mbox {\rm for every}}\hspace{5mm}{g^\prime}\in F_1 .
$$
Which along with (\ref{e2.36}) and  (\ref{e2.38}) implies
$$
\vert \Vert  b_x \eta  \Vert^p - \Vert   b \eta  \Vert^p    \vert < \delta \vert F_1 \vert
$$
and this  (in view of the arbitrariness of $\varepsilon$ and
$\delta$ and recalling  (\ref{e2.35}) and  (\ref{e2.30})) finishes
the verification of  (\ref{e2.28})
\mbox{}\qed\mbox{}\\

Before proceeding to the  next step in the consideration of our example
we have to recall an important for our further purposes notion of amenability.
\begin{I}
\label{2.8a}
{\bf Amenable group.}
\em
Let $G$ be a certain group.
A state $m$ on $l^\infty (G)$ (that is a linear positive functional having the norm equal to 1)
is called a {\em left invariant mean} if
$$
m(f(s^{-1}g)) = m(f(g))
$$
for every $f\in l^\infty (G),$ $s\in G$.\\
A group $G$ is called {\em amenable} if there exists a left invariant mean on  $l^\infty (G)$.\\

In the notion of amenability just defined no topological structure
on the group was presumed thus
$G$ was considered as a {\em discrete} group.\\
The general definition of amenability for locally compact
 topological groups as well as a lot of
information on amenability one can find for example in \cite{Grleaf}.

The definition of amenability presented above used a {\em left}
invariant mean and therefore it is reasonable to consider the
corresponding group as being {\em left} amenable. If one uses {\em
right} invariant means (that is functionals satisfying the
condition
$$
m(f(gs)) = m(f(g)) \hspace{2mm})
$$
then   the corresponding notion of amenability is natural to call the {\em right} amenability.\\

In fact both these notions coincide that is (see \cite{Grleaf}, Lemma 1.1.1.):\\

{\em there exists a left invariant mean on $l^\infty (G)$ iff there exists a
right invariant mean.}

Many  frequently  used groups are amenable, for example\\
\begin{quote}
every abelian group is amenable,\\

every finite group is amenable,\\

every subgroup of an amenable group is amenable,\\

if $S$ is a normal subgroup of $G$ such that $H$ and $G\slash H$ are amenable then
$G$ is amenable (so in particular every solvable group is amenable).\\
\end{quote}
\end{I}
\begin{I}
\label{2.8b}
{\bf Folner's condition.}
\em
There are a number of different criteria for a group $G$ to be amenable.
We shall use the following {\em Folner's condition} (see \cite{Grleaf}, Theorem 3.6.1.):\

{\em a group $G$ is amenable iff \\

for every finite set $K\subset G$ and $\varepsilon > 0$ there exists a
non-empty finite set $U$ such that
$$
\frac {\vert (U\cdot s) \triangle U \vert}{U}  <
 \varepsilon \hspace{5mm}{\mbox  for\hspace{2mm} every}\hspace{5mm} s\in K
$$
where $\vert U \vert $ is the number of elements in $U$.
}
\end{I}
Now the next step in the consideration of our example is
\begin{La}
\label{2.8}
Let $B(A, T_g)$ be the algebra described in \ref{2.6}
and $B({\bar A}, V_g )$ be the corresponding regular representation for $A$
and $\{ {\hat T}_g \}_{g \in G}$ in the space
$H = l^p (G, L^{p}_\mu  (X, E))$. If $G$ is amenable then the mapping
$$
B({\bar A}, V_g ) \to B(A, T_g)
$$
 generated by the mappings
$$
{\bar a}\to a, \hspace{5mm}a\in A
$$
$$
V_g  \to T_g ,  \hspace{5mm}g\in G
$$
is norm decreasing.
\end{La}
{\bf Proof:} Consider the operators
$$
b =  \sum_{g\in F} a_g T_g , \hspace{5mm} \vert F \vert < \infty
$$
and
$$
{\bar b} =  \sum_{g\in F} {\bar a}_g V_g .
$$
We have to show that
\begin{equation}
\label{e2.41}
\Vert {\bar b}\Vert \ge \Vert b \Vert
\end{equation}
Take any $\varepsilon > 0$ and let $f \in  L^{p}_\mu  (X, E)$ be a certain vector such that
\begin{equation}
\label{e2.42}
\Vert f \Vert = 1
\end{equation}
and
\begin{equation}
\label{e2.43}
\Vert bf \Vert > \Vert b \Vert   - \varepsilon
\end{equation}
Choose any finite set $M\subset G$ and define the vector $\eta^{M} \in H$ in the following way
\begin{equation}
\label{e2.44}
\eta^{M}  (g) =
\left \{ \begin{array}{ll}
T_g f, \hspace{5mm} g\in M \\
0, \hspace{5mm} g\notin M\\
\end{array} \right.
\end{equation}
Observe that if for some $h\in G$ we have
$$
\left[  \bigcup_{g\in F}h\cdot g \right]\subset M
$$
then
$$
\Vert (\sum_{g\in F}{\bar a}_g V_g \eta^{M} )(h) \Vert =
\Vert \sum_{g\in F}{\hat T}_h (a_g ) \eta^{M} (hg) \Vert =
$$
\begin{equation}
\label{e2.45}
\Vert \sum_{g\in F} T_h a_g T^{-1}_{hg}f  \Vert =
\Vert \sum_{g\in F}a_g T_g f  \Vert
\end{equation}
Now recalling the  Folner's condition \ref{2.8b}  choose
 finite sets $M_n ,\hspace{2mm}n = 1,2,...$ with the property
$$
\frac{\vert (M_n \cdot g)\triangle M_n \vert}{\vert  M_n \vert}\longrightarrow_{n\to\infty}0,
\hspace{5mm}
g\in F
$$
which implies
\begin{equation}
\label{e2.46}
\frac{\vert \cap_{g\in F}M_n \cdot g \vert}{\vert  M_n \vert}\longrightarrow_{n\to\infty}1
\end{equation}
and consider the vectors $\eta^{M_n}$ given by (\ref{e2.44}).\\
The property (\ref{e2.46}) alon with (\ref{e2.45}) means that
$$
\frac{\Vert  {\bar b}(\eta^{M_n}) \Vert}{\Vert   \eta^{M_n} \Vert }
 \longrightarrow_{n\to\infty}\Vert bf \Vert
$$
And this in view of (\ref{e2.43}) and the arbitrariness of
$\varepsilon$ proves (\ref{e2.41}).
\mbox{}\qed\mbox{}\\

The next lemma closes the investigation circle of the example
considered.
\begin{La}
\label{2.9}
For any element
\begin{equation}
\label{e2.47}
{\bar b} = \sum_{g\in F}{\bar a}_g V_g \in B({\bar A}, V_g ) ,
\hspace{5mm} \vert F \vert < \infty
\end{equation}
we have
\begin{equation}
\label{e2.48}
\Vert {\bar b}\Vert = \sup_{X} \Vert b_x \Vert
\end{equation}
where $b_x , x\in X$ is defined by (\ref{e2.24})- (\ref{e2.26}).
\end{La}
{\bf Proof:} Fix any $\varepsilon >0$ and let $\eta = \{ \eta_g \}_{g \in G} \in H$ be
a certain vector with the properties
\begin{equation}
\label{e2.49}
\Vert \eta \Vert = 1
\end{equation}
and
\begin{equation}
\label{e2.50}
\Vert {\bar b}\eta  \Vert > \Vert {\bar b} \Vert   - \varepsilon
\end{equation}
Without the loss of generality we can assume that there exists a finite
set $M\subset G$ such that
$$
\eta_g = 0 \hspace{5mm}{\rm when} \hspace{5mm}g\notin M .
$$
Recalling the natural isomorphism
\begin{equation}
\label{e2.51}
l^p (G, L^{p}_\mu  (X, E)) \cong L^{p}_{\mu \otimes \mu^\prime}(X\times G, E)
\cong L^{p}_\mu (X,l^p (G,E))
\end{equation}
(here $ \mu^\prime$ is the discrete measure on $G$)\\
we identify the vector $\eta =  \{ \eta_g \} \in l^p (G, L^{p}_\mu  (X, E)) , \hspace{2mm}
\eta_g \in L^{p}_\mu  (X, E)$ with the vector-function
 $ \{ \eta (x) \}_{x\in X}\in L^{p}_\mu (X,l^p (G,E))$.
With this identification we have
\begin{equation}
\label{e2.52}
({\bar b}\eta )(x) = b_x \eta (x)
\end{equation}
and it follows that
$$
\Vert {\bar b}\eta  \Vert \le \sup_X \Vert b_x \Vert \cdot \Vert \eta \Vert =
\sup_X \Vert b_x \Vert
$$
which in view of (\ref{e2.50}) and the arbitrariness of
$\varepsilon$  proves the inequality
$$
\Vert {\bar b}  \Vert \le \sup_X \Vert b_x \Vert .
$$
The opposite inequality can be established by the obvious modification of the argument
 used in the
proof of Lemma \ref{2.7}.
Namely the algebra $\bar A$ can be considered as a subalgebra of the algebra
${\sf A} = C(X\times G, L(E))$ and the action of $G$ on $X\times G$ induced by the automorphisms
$V_g (\cdot ) V^{-1}_g$
$$
{\bar t}_g (x, \tau ) = (x, \tau g^{-1}), \hspace{5mm} (x,\tau )\in X\times G
$$
is evidently topologically free.\\
Thus Lemma \ref{2.9} is proved.
\mbox{}\qed\mbox{}\\
We can summarize the results obtained in
\begin{Tm}
\label{2.10}
Let $B(A, T_g)$ be the algebra described in \ref{2.6}
and $B({\bar A}, V_g )$ be the corresponding regular representation for $A$
and $\{ {\hat T}_g \}_{g \in G}$ in the space
$H = l^p (G, L^{p}_\mu  (X, E))$.\\
 If $G$ is amenable and acts topologically freely then
$$
 B(A, T_g)\cong B({\bar A}, V_g )
$$
where the isomorphism is given by
$$
b = \sum_{g\in F}a_g T_g \longleftrightarrow   \sum_{g\in F}{\bar a}_g V_g = {\bar b},
\hspace{5mm} \vert F \vert < \infty
$$
and in particular the algebra $ B(A, T_g)$ possesses the properties (*) and (**) and
(\ref{e1.36}).
\end{Tm}
{\bf Proof:} Just apply Lemmas \ref{2.7} - \ref{2.9} and the results of Example \ref{1.17} .
\mbox{}\qed\mbox{}\\
\section{ Example 3. Operators in $L^{\infty}_\mu   (X,E)$}
\setcounter{equation}{0}
\begin{I}
\label{2.16} \em Let $X$ be a completely regular space and $\mu$
be a certain Borel measure on $X$ with the support equal to the
whole of $X$. Consider the Banach space
$D= L^{\infty}_\mu   (X,E)$ where $E$ is a Banach space.\\
Let $A= C(X, L(E))$ be the algebra of multiplication operators defined by (\ref{e2.1})
and  $\{ t_g \}_{g \in G}$  be a group of homeomorphisms of $X$ preserving
 the equivalence class of $\mu$.
We denote by $T_g$ the isometry of $D$ given by
\begin{equation}
\label{e2.70}
(T_g f)(x) = f(t^{-1}_g (x))
\end{equation}
Consider the algebra $B(A, T_g )\subset L(D)$ generated by $A$ and $\{  T_g \}_{g \in G}$.
\end{I}
In fact this algebra prove to have the same properties as those obtained when considering
Example 1.
\begin{Tm}
\label{2.17}
Let $B(A, T_g )$ be the algebra introduced above.
If $G$ acts topologically freely then the norm
$\Vert  \sum_{g\in F}a_g T_g \Vert  ,  \hspace{2mm} (\vert F \vert < \infty )$
is calculated via (\ref{e2.5})
\end{Tm}
{\bf Proof:} Since $C(X,E) \subset L^{\infty}_\mu   (X,E)$ it follows
(by Theorem \ref{2.2}) that
$$
\Vert  \sum_{g\in F}a_g T_g \Vert_{L^{\infty}_\mu   (X,E)} \ge
{\sup_X} \sup_{  \{  f_g  \}_{g\in F} \in B_F (E) }
\Vert \sum_{g\in F} a_g  (x) f_g  \Vert
$$
To obtain the opposite inequality just note that for any $f\in L^{\infty}_\mu   (X,E)$ with
$\Vert f \Vert =1$ (which means that without the loss of generality it can be assumed that
$\Vert f \Vert \le1$ for every $x\in X$) we have
$$
\Vert \sum_{g\in F}a_g T_g f \Vert  = {\rm esssup}_X \Vert \sum_{g\in F} a_g  (x) f (t^{-1}_g (x))  \Vert
\le
$$
$$
\sup_X \Vert \sum_{g\in F} a_g  (x) f (t^{-1}_g (x))  \Vert \le
{\sup_X} \sup_{  \{  f_g  \}_{g\in F} \in B_F (E) }
\Vert \sum_{g\in F} a_g  (x) f_g  \Vert
$$
And the proof is complete.
\mbox{}\qed\mbox{}\\
\begin{Rk}
\label{2.18}
\em (cf. Remark \ref{2.3} (1)).
If $E ={\mathbb C}$ and $G$ acts topologically freely then we have
$$
\Vert \sum_{g\in F}a_g T_g  \Vert = \sup_X \vert \sum_{g\in F} a_g  (x) \vert
$$
\end{Rk}
Now the analogue to Lemma \ref{2.4} and Corollary \ref{2.5} is
\begin{La}
\label{2.19}
Let  $B(A, T_g )$ be the algebra described in \ref{2.16} and $B({\bar A}, V_g )$ be
the corresponding regular representation in the space $H = l_0 (G, L^{\infty}_\mu   (X,E) )$
(or $l^\infty  (G, L^{\infty}\mu   (X,E) )$).\\
If $G$ acts topologically freely then
$$
B(A, T_g )\cong B({\bar A}, V_g )
$$
where the isomorphism is generated by the mappings
$$
a \to {\bar a}, \hspace{5mm} a\in A
$$
$$
T_g \to V_g ,  \hspace{5mm}g\in G
$$
and in particular $B(A, T_g )$ possesses the properties (*), (**) and (\ref{e1.36}).
\end{La}
{\bf Proof:} See Lemma \ref{2.4} and Corollary \ref{2.5}.
\mbox{}\qed\mbox{}\\
\section{ Example 4. Operators in $L^{1}_\mu   (X,E)$}
\setcounter{equation}{0}

Our final example deals  with the space $L^1$ and in fact we shall reduce this case to the
previously examined $L^\infty$ situation.\\
We begin with the $L^1 -$analogue to Example 2.
\begin{I}
\label{2.24}
\em
Let $X$ be the space with a measure $\mu$ considered in \ref{2.6},
$D = L^{1}_\mu   (X,E)$ and $A = C(X, L(E))$ be the algebra of operators defined by (\ref{e2.1})
and $\{ T_g \}_{g\in G}$ be defined by (\ref{e2.23}) with $p=1$.\\
Let  $B(A, T_g )\subset L(D)$ be the algebra generated by $A$ and  $\{ T_g \}_{g\in G}$.\\
The space $L^{\infty}_\mu   (X,E^\ast )$ ($E^\ast$ is the dual
space of $E$) is isometrically imbedded in $[L^{1}_\mu   (X,E)
]^\ast$ and for every $\xi \in L^{\infty}_\mu   (X,E^\ast )$ the
value of the corresponding linear functional on $f\in L^{1}_\mu
(X,E) $ is equal to
$$
<f, \xi > = \int_X <f(x), \xi (x)> d\mu
$$
{\em Remark.} It is worth mentioning that in general
$$
[L^{1}_\mu   (X,E)  ]^\ast  \neq L^{\infty}_\mu   (X,E^\ast )
$$
(see for example \cite{Bourb}, Ch. VI, \S 2, 6).\\
The equality here is valid for example when $\mu$ is a discrete
measure or when $E$ is finite dimensional (for some other cases see
 \cite{Bourb}, Ch. VI, \S 2, Exercise 21).\\
But in any case for the spaces considered the following equalities (\ref{e2.79})
and (\ref{e2.80}) are obviously true
\begin{equation}
\label{e2.79}
\Vert f \Vert_{L^1} = \sup_{\Vert \xi \Vert_{L^\infty} = 1} \vert  <f,\xi > \vert
\end{equation}
and
\begin{equation}
\label{e2.80}
\Vert \xi \Vert_{L^\infty} = \sup_{\Vert f \Vert_{L^1} = 1} \vert  <f,\xi > \vert
\end{equation}
\end{I}
\begin{I}
\label{2.25}
\em
For any operator $a\in A$ let us define the operator $a^\natural \in L(L^{\infty}_\mu
 (X,E^\ast ) )$
given by
\begin{equation}
\label{e2.81}
(a^\natural \xi ) (x) = [a(x)]^\ast \xi (x)
\end{equation}
We emphasize once more that while $ [a(x)]^\ast $ is the adjoint to $a(x)$ the operator
$a^\natural $ is in general {\em not necessarily} adjoint to $a$ and we shall call
$a^\natural $ the  {\em formally adjoint} to $a$.\\
For any $T_g$ defined by (\ref{e2.23}) (with $p=1$) we define the operator
$T^{\natural}_g \in L (L^{\infty}_\mu   (X,E^\ast ) )$ by the formula
\begin{equation}
\label{e2.82}
(T^{\natural}_g \xi ) (x) =  \xi (t_g (x))
\end{equation}
And for the same reasons as before we call this operator the {\em formally adjoint} to $T_g$.\\
Finally if $b\in B(A, T_g )$ has the form
$$
b = \sum_{g\in F}a_g T_g , \hspace{5mm} \vert F \vert < \infty
$$
then the {\em formally adjoint} to $b$ is
\begin{equation}
\label{e2.83}
b^{\natural} = \sum_F T^{\natural}_g a^\natural  =
 \sum_F [{\hat T}^{-1}_g a]^{\natural} T^{\natural}_g
\end{equation}
where
\begin{equation}
\label{e2.84}
([{\hat T}^{-1}_g a]^{\natural} \xi )(x) = [a(t_g (x)) ]^{\ast} \xi (x)
\end{equation}
The useful for our subsequent purposes observation is that for any $f \in L^{1}_\mu   (X,E)$
and
$\xi \in L^{\infty}_\mu   (X,E^\ast )$ we have
\begin{equation}
\label{e2.85}
<bf,\xi  > = <f, b^{\natural} \xi >
\end{equation}
which follows obviously from the definition of $b^\natural$.\\
And (\ref{e2.84}) implies in turn (in view of (\ref{e2.79}) and  (\ref{e2.80}))
$$
\Vert b\Vert_{L^1} = \sup_{\Vert f \Vert_{L^1} =1}
 \sup_{\Vert \xi \Vert_{L^\infty} =1}\vert <bf, \xi > \vert =
$$
\begin{equation}
\label{e2.86}
 \sup_{\Vert f \Vert_{L^1} =1}  \sup_{\Vert \xi \Vert_{L^\infty} =1}
 \vert <f, b^{\natural}\xi  > \vert =
\Vert b^{\natural} \Vert_{L^\infty}
\end{equation}
\end{I}
\begin{I}
\label{2.26}
\em
Now let $\{ t_g \}_{g\in G}$ act topologically freely on $X$. The explicit form of the
formally adjoint operator $(T^{\natural}_g$ (see (\ref{e2.82})) shows that the operators
 $(T^{\natural}_g$  induce the automorphisms acting topologically freely on $A^\natural$
(the algebra of operators $a^\natural$ formally adjoint to the operators $a\in A$).\\
In view of this (\ref{e2.86}),  (\ref{e2.83}) and (\ref{e2.84}) and Theorem \ref{2.17}
imply the following result:\\

{\em
 Let $B(A, T_g )$ be the algebra introduced in \ref{2.24}.
If $G$ acts topologically freely then }
$$
\Vert \sum_{g\in F}a_g T_g  \Vert =  {\sup_X} \sup_{  \{  f_g  \}_{g\in F} \in S_F (E^{\ast}) }
\Vert \sum_{g\in F} [a_g  (t_g (x)]^{\ast} f_g  \Vert =
$$
\begin{equation}
\label{e2.87}
 {\sup_X} \sup_{  \{  f_g  \}_{g\in F} \in B_F (E^{\ast}) }
\Vert \sum_{g\in F} [a_g (t_g (x)]^{\ast}  f_g  \Vert
\end{equation}
\end{I}
\begin{Rk}
\label{2.27}
\em (cf. Remark \ref{2.18}).
If $E={\mathbb C}$ and $G$ acts topologically freely then
$$
\Vert \sum_{g\in F}a_g T_g  \Vert =  {\sup_X}\sum_{g\in F}\vert a_g (t_g (x) \vert
$$
\end{Rk}
Now the analogue to Lemma \ref{2.19} is
\begin{La}
\label{2.28}
Let  $B(A, T_g )$ be the algebra described in \ref{2.24} and  $B({\bar A}, V_g )$ be the
corresponding regular representation in the space $ H = l^1 (G, L^{1}_\mu (X, E))$.
If $G$ acts topologically freely then
$$
B(A, T_g ) \cong B({\bar A}, V_g )
$$
where the isomorphism is generated by the mappings
$$
a \to {\bar a}, \hspace{5mm} a\in A
$$
$$
T_g \to V_g ,  \hspace{5mm} g\in G
$$
and in particular $B(A, T_g )$ possesses the properties (*) and (**) and (\ref{e1.36}).
\end{La}
{\bf Proof:} The idea here is the same as in the proof of Lemma \ref{2.4} and Corollary \ref{2.5}.
The only difference is the usage of (\ref{e2.87}) instead of (\ref{e2.5}).
\mbox{}\qed\mbox{}\\
Now we would like to observe certain interrelations between the examples considered and the
Isomorphism Theorem (\cite{AnLeb}, Corollary 12.17).
 \begin{Rk}
\label{2.33}
\em
(1) Observe that in Examples 1, 3, 4 we did not use any information about
the group $G$ thus the group in these examples is
{\em not necessarily amenable}.\\

(2) The essentially different picture is drawn in Example 2.\\
Here\\
(i) if $G$ acts topologically freely then \\
$B({\bar A}, V_g )$ is a representation of $B(A, T_g )$ for {\em any} $G$
({\em not necessarily amenable})\\
 (Lemma \ref{2.7} and Lemma \ref{2.9}).\\

While\\
(ii) if $G$ is amenable then\\
 $B(A, T_g )$ is a representation of $B({\bar A}, V_g )$ for an {\em arbitrary} action of $G$
({\em not necessarily topologically free}).\\
(Lemma \ref{2.8}).\\

Thus in these examples the topological freedom of the action of $G$ and the amenability
of $G$ are lying in a sense opposite each other.\\
If $G$ acts topologically freely then  $B(A, T_g )$  is 'larger' than $B({\bar A}, V_g )$
(see (i)).\\
And\\
if $G$ is amenable then  $B({\bar A}, V_g )$ is 'larger' than  $B(A, T_g )$ (see (ii)).\\

Both these algebras 'coincide' if $G$ acts topologically freely
and is amenable
(Theorem \ref{2.10}).\\

(3) Consideration of Example 2 leads to certain {\em Isomorphism Theorems} which
(just as it was done in  \cite{AnLeb}, Corollary 12.17) establish
the isomorphism between {\em essentially
spatially different} operator algebras (thus wiping off the spaces where these operators  act).\\
For example.\\
Let $\mu_1$ and $\mu_2$ be two Borel measures on a completely
regular space $X$ each having the support equal to the whole of
$X$ and let $\{ t_g \}_{g\in G}$ be a group of homeomorphisms of
$X$ preserving the equivalence classes of $\mu_1$ and $\mu_2$.
Consider the spaces $D_i = L^{p}_{\mu_i} (X, E),\hspace{2mm}i=1,2
$. And Let $A_i = C(X, L(E) \subset  L(D_i )), \hspace{2mm}i=1,2$
be the algebras of multiplication operators defined by
(\ref{e2.1}) and $T^{i}_g ,\hspace{2mm}i=1,2$ be the isometries of
$D_i$ defined by (\ref{e2.23}) (with $\mu =\mu_i$) and $B(A_i
,T^{i}_g ), \hspace{2mm}i=1,2$ be the algebras generated by $A_i$
and $\{ T^{i}_g \}_{g\in G}$.\\
The Isomorphism Theorem we have in mind is stated as follows:\\

{\em
If $G$ is amenable and acts topologically freely then $B(A_1 ,T^{1}_g )$
and $B(A_2 ,T^{2}_g )$ are isomorphic  (as Banach algebras)
and the isomorphism is established by the natural isomorphism
$$
A_1 \cong A_2
$$
and the mapping
$$
T^{1}_g \to T^{2}_g
$$
}
{\bf Proof:} In view of Lemmas \ref{2.7}, \ref{2.8} and \ref{2.9} we have
\begin{equation}
\label{e2.89}
B(A_i ,T^{i}_g ) \cong { \bigoplus_{x\in X}} \pi_x ( B(A_i ,T^{i}_g ) )
\end{equation}
where $\pi_x$ is defined by (\ref{e2.24}) - (\ref{e2.26}).\\
Now observe that by the explicit form of $\pi_x$ the right-hand part in (\ref{e2.89})
does not depend on $\mu_i , \hspace{2mm}i=1,2$ but only on $\{ t_g \}_{g\in G}$
which proves the statement.
\mbox{}\qed\mbox{}\\
\end{Rk}
\section{Interpolation}
\setcounter{equation}{0}

We have not obtained the explicit formula for the norm of operator
$\sum a_g T_g$ in Example 2. But in fact the formulae proved while considering the other examples
and the Riesz-Thorin interpolation theorem makes it possible to write out the useful
estimates for the norm of $\sum a_g T_g$ in all the cases.\\

To present  the result we have in mind let us introduce first the space
$L^{0}_{\mu} (\Omega , E)$ which is the closed subspace of $L^{\infty}_{\mu} (\Omega , E)$
generated by the functions having the supports of finite measure that is the functions $f$
such that
$$
 f\in L^{\infty}_{\mu} (\Omega , E)
$$
and
$$
\mu (\{ x: f(x) \neq 0 \}) < \infty
$$
Obviously if $\mu$ is a discrete measure then
$$
L^{0}_{\mu} (\Omega , E) = l_0 (\Omega , E)
$$
The next Theorem \ref{2.35} is in fact a certain vector valued variant of the classical
{\em Riesz-Thorin interpolation theorem} (see \cite{BerLof}, 1.1.1)
\begin{Tm}
\label{2.35}
Let $V$ be a linear operator belonging to $L(D_1 ) $ and  $L(D_2 ) $ where
$D_i = L^{p_i}_{\mu} (\Omega , E)$ and $p_1 \in [1, \infty  )$ while
$p_2 \in [1, \infty  ) \cup \{ 0\}$.
For any number $\Theta , \hspace{2mm} 0< \Theta  < 1$  let $p$ be the number defined by
$$
{1\over p} = \frac{1 - \Theta}{p_1} +  \frac{\Theta}{p_2}\hspace{5mm}{\rm when}
\hspace{2mm} p_2 \neq 0
$$
or
$$
{1\over p} = \frac{1 - \Theta}{p_1} \hspace{5mm}{\rm when}
\hspace{2mm} p_2 = 0
$$
and let $D = L^{p}_{\mu} (\Omega , E)$.\\
Then the operator $V$ belongs to $L(D)$ and
\begin{equation}
\label{2.90}
\Vert V\Vert_D \le \Vert V\Vert^{1-\Theta}_{D_1} \cdot \Vert V\Vert^{\Theta}_{D_2}
\end{equation}
\end{Tm}
{\bf Proof:}  See \cite{BerLof}, 5.1.2 and 4.1.2 .
\mbox{}\qed\mbox{}\\
Now let us apply this result to the operators considered in Example 2.\\

Let $b \in B(A, T_g ) \subset L(L^{1}_\mu   (X,E)), \hspace{2mm} 1\le p\le \infty $
be an operator of the form
\begin{equation}
\label{2.91}
b = \sum_{g\in F}a_g T_g , \hspace{5mm} \vert F \vert < \infty
\end{equation}
where for $1\le p < \infty $ $B(A, T_g )$ is that described in \ref{2.6}
and for $p = \infty$ in \ref{2.16}.\\
Let also ${\bar b} \in B({\bar A}, V_g )$ be the operator of the form
\begin{equation}
\label{2.92}
{\bar b} = \sum_{g\in F}{\bar a}_g V_g
\end{equation}
where $B({\bar A}, V_g )$ is the corresponding regular representation for
$A$ and $\{ {\hat T}_g \}_{g\in G}$ respectively in the space
$l^p (G, L^{p}_\mu   (X,E) ), \hspace{2mm}1\le p\le \infty$
(in the case $p=\infty$ we can use the space $l^p (G, L^{p}_\mu   (X,E) )$ as well).\\
We shall denote by $\Vert b \Vert_p$ the norm of the operator $b$ as acting
 in the space $L^{p}_\mu   (X,E)$. And in the analogous way we shall introduce the notation
 $\Vert{\bar b} \Vert_p$.\\
The important observation here is that the operator $\bar b$ is given by the same formula
for all the cases $1\le p \le\infty$ (this will give us
the possibility to apply Theorem \ref{2.35}).\\
Recall that the results of Examples 2 - 4 show that if $G$ is amenable and acts topologically
freely then for every $p\in [1,\infty ]$
$$
\Vert b \Vert_p  = \Vert {\bar b} \Vert_p
$$
(in the cases $p=1$ or $\infty$ the amenability of $G$ is not necessary).\\
This equality along with the observation on the structure of $b$ just mentioned
and Theorem \ref{2.35}
imply the following result.
\begin{La}
\label{2.37}
If $G$ is amenable and acts topologically freely on $X$ then for any $p\in (1,\infty )$ we have
\begin{equation}
\label{e2.93}
\Vert b \Vert_p \le \Vert b \Vert^{1\over p}_1 \cdot \Vert b \Vert^{1 - {1\over p}}_\infty
\end{equation}
where $\Vert b \Vert_1$ is given by (\ref{e2.87}) and  $\Vert b \Vert_\infty$ by (\ref{e2.5}).
\end{La}
\begin{I}
\label{2.38}
\em
Looking through Example 2 a bit more thoroughly we can strengthen the estimate
(\ref{e2.93}). Namely recall that if $G$ is amenable and acts topologically freely the
rsults of Examples 2, 3 and 4  also show that
\begin{equation}
\label{e2.94}
\Vert b \Vert_p = \Vert {\bar b} \Vert_p = \sup_X \Vert b_x \Vert_p
\end{equation}
where $b_x$ is defined by (\ref{e2.24}) - (\ref{e2.26})
(and we consider this operator as acting in the corresponding space $l^p (G,E), \hspace{2mm}
1\le p \le \infty$ or in the case $p=\infty $ we also can use the space $l_0 (G,E)$).\\
For the operators $b_x$ Theorem \ref{2.35} imlies\\
{\em
for every $p\in (1,\infty )$
\begin{equation}
\label{2.95}
\Vert b_x \Vert_p \le \Vert b_x \Vert^{1\over p}_1 \cdot \Vert b_x \Vert^{1 - {1\over p}}_\infty =
\Vert b_x \Vert^{1\over p}_1 \cdot \Vert b_x \Vert^{1 - {1\over p}}_0
\end{equation}
(here $ \Vert b_x \Vert_0$ is the norm of $b_x$ in $l_0 (G,E)$).
}\\
Which means in turn in view of (\ref{e2.94}) that\\
{\em
if $G$ is amenable and acts topologically freely on $X$ then for any}
$p\in (1, \infty )$
$$
\Vert b \Vert_p =  \Vert {\bar b} \Vert_p \le \sup_X
\Vert b_x \Vert^{1\over p}_1 \cdot \Vert b_x \Vert^{1 - {1\over p}}_\infty
$$
We also recall that the results of Examples 1 and 4 show that
$$
 \Vert b_x \Vert_\infty =  \Vert b_x \Vert_0 =
\sup_{h\in G}
 \sup_{  \{  f_g  \}_{g\in F} \in S_F (E) }
\Vert \sum_{g\in F}a_g  (t^{-1}_h (x) f_g  \Vert  =
$$
$$
\sup_{h\in G}
 \sup_{  \{  f_g  \}_{g\in F} \in B_F (E) }
\Vert \sum_{g\in F}a_g  (t^{-1}_h (x) f_g  \Vert
$$
and
$$
 \Vert b_x \Vert_1 =
\sup_{h\in G}
 \sup_{  \{  f_g  \}_{g\in F} \in S_F (E^\ast ) }
\Vert \sum_{g\in F}[a_g  (t_{g\cdot h^{-1}} (x)]^{\ast} f_g  \Vert  =
$$
$$
\sup_{h\in G}
 \sup_{  \{  f_g  \}_{g\in F} \in B_F (E^\ast ) }
\Vert \sum_{g\in F}[a_g  (t_{g\cdot h^{-1}} (x)]^{\ast} f_g  \Vert  .
$$
\end{I}

\end{document}